\documentclass[11pt]{article}

\usepackage{amsfonts}
\usepackage{amssymb}
\usepackage{enumerate}
\newtheorem{st}      {Theorem}
\newtheorem{prop}  {Proposition}
\newtheorem{lem} {Lemma}
\newtheorem{cor} {Corollary}

\newtheorem{exam}{Example}
\newtheorem{num}  {}

\newcommand{\rmap}{\longrightarrow}
\newcommand{\Boxe}{\raisebox{.8ex}{\framebox}}

\newcommand{\xx}{\ensuremath{\mathcal{X}}}

\newcommand{\A}{\ensuremath{\mathcal{A}}}

\newcommand{\F}{\ensuremath{\mathcal{F}}}

\newcommand{\ps}{{\raise 1pt\hbox{\tiny (}}}

\newcommand{\pss}{{\raise 1pt\hbox{\tiny [}}}
\newcommand{\pdd}{{\raise 1pt\hbox{\tiny ]}}}
\newcommand{\pd}{{\raise 1pt\hbox{\tiny )}}}

\newcommand{\bs}{{\raise 1pt\hbox{\tiny [}}}
\newcommand{\bd}{{\raise 1pt\hbox{\tiny ]}}}

\def\cross{\mathinner{\mathrel{\raise0.8pt\hbox{$\scriptstyle>$}}
                 \joinrel\mathrel\triangleleft}}

\def\compose{{\raise 1pt\hbox{$\scriptscriptstyle\circ$}}}
\def\dcross{{\raise 0.5pt\hbox{$\scriptscriptstyle\boxtime$}}}

\advance\textheight by \topskip
\usepackage[all]{xy}
\usepackage{epsf}
\usepackage{amsfonts}
\usepackage{amssymb}
\usepackage{amsmath}

\begin{document}

\title{Connections up to homotopy and characteristic classes  \thanks{Research supported by NWO}}
\author {Marius Crainic}
\date {Department of Mathematics, Utrecht University, The Netherlands}
\pagestyle{myheadings}
\maketitle


\section*{Introduction}

\hspace*{.3in}The aim of this note is to clarify the relevance of ``connections up to homotopy'' \cite{Crai, ELW} to the theory of characteristic classes,
and to present an application to the characteristic classes of Lie algebroids \cite{Cra, ELW, Fer1} (and of Poisson manifolds
in particular \cite{Fer2, We}).\\
\hspace*{.3in}We have already remarked \cite{Crai} that such connections up to homotopy can be used to compute the classical Chern characters.
Here we present a slightly different argument for this, and then proceed with the discussion of the flat characteristic classes.  
In contrast with \cite{Crai}, we do not only recover the classical characteristic classes
(of flat vector bundles), but we also obtain new ones. The reason for this is that ($\mathbb{Z}_2$-graded) non-flat vector 
bundles may have flat connections up to homotopy. As we shall explain here, in this category fall e.g. the 
characteristic classes of Poisson manifolds \cite{Fer2, We}.\\
\hspace*{.3in}As already mentioned in \cite{Crai}, one of our motivations is to understand the intrinsic characteristic classes
for Poisson manifolds (and Lie algebroids) of \cite{Fer1, Fer2}, and the connection with the characteristic classes of representations \cite{Cra}.
Conjecturally, Fernandes' intrinsic characteristic classes \cite{Fer1} are the characteristic classes \cite{Cra} of the ``adjoint
representation''. The problem is that the adjoint representation is a ``representation up to homotopy'' only. Applied to Lie algebroids,
our construction immediately solves this problem: it extends the characteristic classes of \cite{Cra} from representations
to representations up to homotopy, and shows that the intrinsic characteristic classes \cite{Fer1, Fer2} are indeed the ones associated to the adjoint representation \cite{ELW}.\\
\hspace*{.3in}I would like to thank J. Stasheff and A. Weinstein for their comments on a preliminary version of this paper.

\section*{Non-linear connections}


\hspace*{.3in}Here we recall some well-known properties of connections
on vector bundles. Up to a very slight novelty (we allow non-linear connections),
this section is standard \cite{Chern} and serves to fix the notations.\\ 
\hspace*{.3in}Let $M$ be a manifold, and let $E= E^0\oplus E^1$ be a super-vector bundle over $M$. We now consider $\mathbb{R}$-linear 
operators 
\begin{equation}\label{connection}
 \xx(M)\otimes \Gamma E \rmap \Gamma E, \ \ (X, s)\mapsto \nabla_{X}(s)
\end{equation}
which satisfy 
\[ \nabla_{X}(fs)= f\nabla_{X}(s) + X(f) s \]
for all $X\in \xx(M)$, $s\in \Gamma E$, and $f\in C^{\infty}(M)$, and which preserve the grading of $E$. 
We say that $\nabla$ is a {\it non-linear connection} if $\nabla_{X}(V)$ is local in $X$. This is just a relaxation
of the $C^{\infty}(M)$-linearity in $X$, when one recovers the standard notion of (linear) connection. 
The curvature $k_{\nabla}$ of a non-linear connection $\nabla$ is defined 
by the standard formula
\begin{equation}\label{curvature}
k_{\nabla}(X, Y)= [\nabla_{X}, \nabla_{Y}]- \nabla_{[X, Y]}: \Gamma E\rmap \Gamma E\ .
\end{equation}
\hspace*{.3in}A {\it non-linear differential form}\footnote{as in the case of connections,
the non-linearity referes to $C^{\infty}(M)$-non-linearity. As pointed out to me, the terminology might be misleading. Betters
names would probably be ``higher order connections'' and ``jet-forms''} on $M$ is an antisymmetric ($\mathbb{R}$-multilinear) map
\begin{equation}\label{omega}
\omega: \underbrace{\xx(M)\times \ldots \times \xx(M)}_{n} \rmap C^{\infty}(M)
\end{equation}
which is local in the $X_{i}$'s.  
It is easy to see (and it has been already remarked in \cite{Crai}) that many of the usual operations 
on differential forms do not use the $C^{\infty}(M)$-linearity, hence they apply to non-linear forms as well. In particular  
we obtain the algebra $(\A_{\text{nl}}(M), d)$ of non-linear forms endowed with De Rham operator. (This defines a contravariant functor
from manifolds to dga's.) Considering $\Gamma E$-valued operators instead, we obtain a version with coefficients, denoted
$\A_{\text{nl}}(M; E)$. Note that a non-linear connection $\nabla$ can be viewed as an operator
$\A_{\text{nl}}^{0}(M; E) \rmap \A_{\text{nl}}^{1}(M; E)$ which has a unique extension to an operator
\[ d_{\nabla}: \A_{\text{nl}}^{*}(M; E)\rmap \A_{\text{nl}}^{*+1}(M; E) \]
satisfying the Leibniz rule. Explicitly,
\begin{eqnarray}\label{differential}
d_{\nabla}(\omega)(X_1, \ldots , X_{n+1}) & = & \sum_{i<j}
(-1)^{i+j}\omega([X_i, X_j], X_1, \ldots , \hat{X_i}, \ldots ,
\hat{X_j}, \ldots X_{n+1})) \nonumber \\
 & + & \sum_{i=1}^{n+1}(-1)^{i+1}
\nabla_{X_i}\omega(X_1, \ldots, \hat{X_i}, \ldots , X_{n+1}) .
\end{eqnarray}
\hspace*{.3in}We now recall the definition of the (non-linear) connection on $\text{End}(E)$ induced by $\nabla$. 
For any $T\in \Gamma \text{End}(E)$, the operators $[\nabla_X, T]$ acting on $\Gamma(E)$ are
$C^{\infty}(M)$-linear, hence define elements $[\nabla_X, T]\in \Gamma \text{End}(E)$. The desired connection is then
$\nabla_{X}(T)= [\nabla_X, T]$. Clearly $k_{\nabla}\in \A_{\text{nl}}^{2}(M; \text{End}(E))$, and one has Bianchi's identity $d_{\nabla}(k_{\nabla})= 0$.\\
\hspace*{.3in}We will use the algebra $\A_{\text{nl}}(M; \text{End}(E))$
and its action on $\A_{\text{nl}}(M; E)$. The product
structure that we consider here is the one which arises from the natural isomorphisms
\[ \A_{\text{nl}}(M; E)\cong \A_{\text{nl}}(M)\otimes_{C^{\infty}(M)}\Gamma(E) \]
and the usual sign conventions for the tensor products (i.e. 
$\omega\otimes x\cdot \eta\otimes y= (-1)^{|x||\eta|} \omega\eta \otimes xy$). 
The usual super-trace on $\text{End}(E)$ induces a super-trace 
\begin{equation}\label{supertrace}
Tr_{s}: (\A_{{\rm nl}}(M;\, {\rm End}(E)), d_{\nabla}) \rmap (\A_{{\rm nl}}(M), d)
\end{equation}
with the property that $Tr_{s}d_{\nabla}= d Tr_{s}$. We conclude (and this is just a non-linear version of the standard
construction of Chern characters \cite{Chern}):

\begin{lem} If $\nabla$ is a non-linear connection on $E$, then
\begin{equation}\label{character}
ch_{p}(\nabla)= Tr_{s}(k_{\nabla}^{p}) \in \A^{2p}_{{\rm nl}}(M) 
\end{equation}
are closed non-linear forms on $M$.
\end{lem}

\hspace*{.1in}Up to a boundary, these classes are independent of $\nabla$. This is an instance of the Chern-Simons construction
that we now recall. Given $k+1$ non-linear connections $\nabla_i$ on $E$ ($0\leq i\leq k$) we form
their affine combination 
$\nabla^{\text{aff}}= (1- t_1- \ldots - t_k)\nabla_0+ t_1\nabla_1+ \ldots + t_{k}\nabla_k$. This is a non-linear connection
on the pullback of $E$ to ${\bf \Delta}^{k}\times M$, where ${\bf \Delta}^{k}= \{ (t_1, \ldots, t_k): t_i\geq 0, \sum t_{i}\leq 1\}$ is the standard $k$-simplex. The classical integration along fibers has a non-linear extension
\begin{equation}\label{integration}
\int_{{\bf \Delta}^{k}}: \A_{\text{nl}}^{*}(M\times {\bf \Delta}^{k}) \rmap \A_{\text{nl}}^{*-k}(M) 
\end{equation}
given by the explicit formula
\[ (\int_{{\bf \Delta}^{k}} \omega)(X_1, \ldots , X_{n-k})= \int_{{\bf \Delta}^{k}}\omega (\frac{\partial}{\partial t_1}, \ldots ,
\frac{\partial}{\partial t_k}, X_1, \ldots , X_{n-k}) dt_1 \ldots dt_{k}\ . \]
We then define
\begin{equation}\label{ChernSimons}
cs_{p}(\nabla_0, \ldots , \nabla_k)= \int_{{\bf \Delta}^{k}} ch_{p}(\nabla^{\text{aff}}) \ .
\end{equation}
Using a version of Stokes' formula \cite{Bott} (or integrating by parts repeatedly) we conclude

\begin{lem} The elements (\ref{ChernSimons}) satisfy
\begin{equation}\label{Stokes} 
dcs_p(\nabla_0, \ldots, \nabla_k)= \sum_{i=0}^{k} (-1)^{i} cs_{p}(\nabla_0, \ldots , \widehat{\nabla_i}, \ldots , \nabla_k)\ .
\end{equation} 
\end{lem}

\section*{Connections up to homotopy and Chern characters}

\hspace*{.3in}From now on, $(E, \partial)$ is a super-complex of vector bundles over the manifold $M$,
\begin{eqnarray}\label{complex}
(E, \partial): \ \  \xymatrix{  E^0\ \ar@<-1ex>[r]_-{\partial} & \ \ E^1 \ar@<-1ex>[l]_-{\partial}\ \ \ \ . }
\end{eqnarray}
We now consider non-linear connections $\nabla$ on $E$ such that $\nabla_X\partial= \partial\nabla_X$ for all $X\in \xx(M)$.
We say that $\nabla$ is a (linear) {\it connection} on $(E, \partial)$ if it also satisfies the identity $\nabla_{fX}(s)= f \nabla_{X}(s)$
for all $X\in \xx(M)$, $f\in C^{\infty}(M)$, $s\in \Gamma E$. The notion of {\it connection up to homotopy} \cite{Crai, ELW} on $(E, \partial)$ is obtained by relaxing the $C^{\infty}(M)$-linearity on $X$ to linearity up to homotopy. In other words we require 
\[ \nabla_{fX}(s)= f \nabla_{X}(s) + [H_{\nabla}(f, X), \partial ]\ , \]
where $H_{\nabla}(f, X)\in \Gamma \text{End}(E)$ are odd elements which are $\mathbb{R}$-linear and local in $X$ and $f$.\\
\hspace*{.3in}We say that two non-linear connections $\nabla$ and $\nabla^{\,'}$ are {\it equivalent}
(or homotopic) if 
\[ \nabla_{X}^{\,'}= \nabla_{X} + [\theta(X), \partial] \]
for all $X\in \xx(M)$, for some $\theta\in \A_{\text{nl}}^{1}(M;\text{End}(E))$ of even degree. We write $\nabla\sim \nabla^{\, '}$.  

\begin{lem}\label{lemmaexist} 
A non-linear connection is a connection up to homotopy if and only if it is equivalent to a (linear) connection.
\end{lem}

{\it Proof: } Assume that $\nabla$ is a connection up to homotopy. Let $U_{a}$ be the domain of local coordinates $x^k$ for $M$, and 
put
\[ \nabla^{a}_{X}= \nabla_{X}+ [u^{a}(X), \partial ]\ ,\]
where $u_{a}\in \A_{\text{nl}}(U_a; \text{End}(E))$ is given by
\[ u_{a}(\sum_{k} f_k \frac{\partial}{\partial x_{k}})= - \sum_{k} H_{\nabla}(f_{k}, \frac{\partial}{\partial x_{k}}) \ ,\]
for all $f_k\in C^{\infty}(U_a)$. Note that $\nabla_{X}$ is linear on $X$. Indeed,
for any two smooth functions $f, g$ and $X= g\frac{\partial}{\partial x_{k}}$ we have
\begin{eqnarray}
 & \nabla_{fX}^{a}- f\nabla_{X}^{a}  = (\nabla_{fX}+ [u^a(fX), \partial ])- f(\nabla_{X} + [u^a(X), \partial ]) =  \nonumber \\
 & = (\nabla_{fg\frac{\partial}{\partial x_{k}}}- [H_{\nabla}(fg, \frac{\partial}{\partial x_{k}}), \partial ])+ f 
(\nabla_{g\frac{\partial}{\partial x_{k}}}- [H_{\nabla}(g \frac{\partial}{\partial x_{k}}), \partial ]) =  \nonumber\\
 & = fg\nabla_{\frac{\partial}{\partial x_{k}}} - f g\nabla_{\frac{\partial}{\partial x_{k}}} =  0 \ . \nonumber
\end{eqnarray}
Next we take $\{ \nu_{a}\}$ to be a partition of unity subordinate to an open cover $\{ U_a\}$ by such coordinate domains
and put $\nabla_{X}^{\, '}= \sum_{a} \nu_{a} \nabla_{X}^{a}$, $u(X)= \sum_{a} \nu_{a} u^{a}(X)$. Then $\nabla^{\, '}= \nabla + [u, \partial]$ is a connection equivalent to $\nabla$. \ \ \ $\Boxe$\\

\begin{lem}\label{lemequiv} If $\nabla_{0}$ and $\nabla_{1}$ are equivalent, then  $ch_{p}(\nabla^{0})= ch_{p}(\nabla^{1})$.
\end{lem}

{\it Proof: } So, let us assume that $\nabla^{1}= \nabla^{0} + [\theta, \partial]$. A simple computation shows that 
\begin{equation}\label{decomp}
k_{\nabla_{1}}= k_{\nabla_{0}}+ [d_{\nabla}(\theta)+ R, \partial] \ ,
\end{equation}
where $R(X, Y)= [\theta(X), [\theta(Y), \partial ]]$. 
We denote by $Z\subset \A_{\text{nl}}(M; \text{End}(E))$ the space of non-linear forms $\omega$ with the property that $[\omega, \partial ]= 0$,
and by $B\subset Z$ the subspace of element of type $[\eta, \partial ]$ for some non-linear form $\eta$. The formula
\[ [\partial ,\omega\eta ]= [\partial , \omega ]\eta + (-1)^{|\omega |}\omega [\partial , \eta ] \]
shows that $ZB\subset B$, hence (\ref{decomp}) implies that $k_{\nabla_{1}}^{p}\equiv k_{\nabla_{0}}^{p}$ modulo $B$. 
The desired equality follows now from the fact that 
$Tr_{s}$ vanishes on $B$.  \ \ \ $\Boxe$\\

\hspace*{.1in}For (linear) connections $\nabla$ on $(E, \partial)$, $ch_{p}(\nabla)$ are clearly (linear) differential forms
on $M$ whose cohomology classes are (up to a constant) the components of the Chern character $Ch(E)= Ch(E^0)- Ch(E^1)$.
Hence an immediate consequence of the previous two lemmas is the following \cite{Crai}

\begin{st}\label{theorem1} If $\nabla$ is a connection up to homotopy
on $(E, \partial)$, then $ch_p(\nabla)= Tr_{s}(k_{\nabla}^{p})$
are closed differential forms on $M$ whose De Rham cohomology classes are (up to a constant) the components of
the Chern character $Ch(E)$.
\end{st}

\section*{Flat characteristic classes} 

\hspace*{.3in}As usual, by flatness we mean the vanishing of the curvature forms. Theorem \ref{theorem1} immediately implies

\begin{cor} 
If $(E, \partial)$ admits a connection up to homotopy which is flat, then $Ch(E)= 0$.
\end{cor}

\hspace*{.1in}As usual, such a vanishing result is at the origin of new ``secondary'' characteristic classes.
Let $\nabla$ be a flat connection up to homotopy. To construct the associated secondary classes we need
a metric $h$ on $E$.  We denote by $\partial^h$ be the adjoint of $\partial$ with respect to $h$.  Using the isomorphism
$E^*\cong E$ induced by $h$ (which is anti-linear if $E$ is complex), $\nabla$ induces an adjoint connection $\nabla^h$ on $(E, \partial^h)$. Explicitly,
\[ L_{X} h(s, t)= h(\nabla_{X}(s), t)+ h( s, \nabla_{X}^{h}(t)) \ .\]
The following describes various possible definitions of the secondary classes, as well as their main properties
(note that the role of $i= \sqrt{-1}$ below is to ensure real classes).


\begin{st}\label{theorem2} Let $\nabla$ be a flat connection up to homotopy on $(E, \partial)$, $p\geq 1$.
\begin{enumerate} [(i)]
\item For any (linear) connection $\nabla_{0}$ on $(E, \partial)$ and any metric $h$, 
\begin{equation}\label{classes} 
i^{p+1}(cs_{p}(\nabla, \nabla_0)+ cs_{p}(\nabla_{0}, \nabla_{0}^{h})+ cs_{p}(\nabla_{0}^{h}, \nabla^{h})) \in \A_{{\rm nl}}^{2p-1}(M)
\end{equation}
are differential forms on $M$ which are real and closed. The induced cohomology classes do not depend on the choice of $h$ or $\nabla_0$,
and are denoted $u_{2p-1}(E, \partial, \nabla)\in H^{2p-1}(M)$.
\item For any connection $\nabla_{0}$ equivalent to $\nabla$, and any metric $h$,
\begin{equation}\label{classes2} 
i^{p+1}cs_{p}(\nabla_{0}, \nabla_{0}^{h}) \in \A^{2p-1}(M)
\end{equation}
are real and closed, and represent $u_{2p-1}(E, \partial, \nabla)$ in cohomology.
\item If $\nabla$ is equivalent to a metric connection (i.e. a connection which is compatible with a metric), then 
all the classes $u_{2p-1}(E, \partial, \nabla)$ vanish.
\item If $\nabla \sim \nabla^{\, '}$, then $u_{2p-1}(E, \partial, \nabla)= u_{2p-1}(E, \partial, \nabla^{\, '})$.
\item If $\nabla$ is a flat connection up to homotopy on both super-complexes $(E, \partial)$ and $(E, \partial^{\,'})$, then $u_{2p-1}(E, \partial, \nabla)=
u_{2p-1}(E, \partial^{\,'}, \nabla)$.
\item Assume that $E$ is real. If $p$ is even then $u_{2p-1}(E, \partial, \nabla)= 0$. If $p$ is odd, then 
for any connection $\nabla_{0}$ equivalent to $\nabla$, and any metric connection $\nabla_{\text{m}}$,
\[ (-1)^{\frac{p+1}{2}}cs_{p}(\nabla_{0}, \nabla_{\text{m}}) \in \A^{2p-1}(M) \]
are closed differential forms whose cohomology classes equal to $\frac{1}{2}u_{2p-1}(E, \partial, \nabla)$.
\end{enumerate}
\end{st}

\hspace*{.1in}Note the compatibility with the classical flat characteristic classes, which 
correspond to the case where $E$ is a graded vector bundle (and $\partial= 0$),
or, more classically, just a vector bundle over $M$. As references for this we
point out \cite{KaTo} (for the approach in terms of frame bundles and Lie algebra cohomology),
and \cite{BiLo} (for an explicit approach which we follow here). For the proof of the theorem we need the following

\begin{lem}\label{lemmacs} Given the non-linear connections $\nabla$, $\nabla_0$, $\nabla_1$,
\begin{enumerate}[(i)]
\item If $\nabla_{0}$ and $\nabla_{1}$ are connections up to homotopy then
$cs_{p}(\nabla_{0}, \nabla_{1})$ are differential forms;
\item If $\nabla_{0}\sim \nabla_{1}$, then $cs_{p}(\nabla_{0}, \nabla_{1})= 0$;
\item For any metric $h$, $ch_{p}(\nabla^h)= (-1)^p \overline{ch_{p}(\nabla)}$
and $cs_{p}(\nabla_{0}^{h}, \nabla_{1}^{h})= (-1)^p\overline{cs_{p}(\nabla_{0}, \nabla_{1})}$.
\end{enumerate}
\end{lem}

{\it Proof:} (i) follows from the fact that Chern characters of connections
up to homotopy are differential forms. For (ii) we use Lemma \ref{lemequiv}. The affine combination $\nabla$
used in the definition of $cs_{p}(\nabla_{0}, \nabla_{1})$ is equivalent to the pull-back $\tilde{\nabla}_{0}$
of $\nabla_{0}$ to $M\times {\bf \Delta}^{1}$ (because $\nabla= \tilde{\nabla}_{0}+ t[\theta, \partial]$),
while $ch_{p}(\tilde{\nabla}_{0})$ is clearly zero. If $h$ is a metric on $E$, a simple computation
shows that $k_{\nabla^h}(X, Y)$ coincides with $-k_{\nabla}(X, Y)^{*}$ where $*$ denotes the adjoint (with respect to $h$).
Then (iii) follows from $Tr(A^{*})= \overline{Tr(A)}$ for any matrix $A$.\ \ \ $\Boxe$ \\

{\it Proof of Theorem \ref{theorem2}:} (i) Let us denote by $u(\nabla, \nabla_0, h)$ the forms (\ref{classes}).
Since $(\nabla_{0}, \nabla^{h}_{0})$ is a pair of connections on $E$, and $(\nabla, \nabla_0)$, 
$(\nabla^h, \nabla^{h}_{0})$ are pairs of connections up to homotopy on $(E, \partial)$ and $(E, \partial^h)$, respectively,
it follows from (i) of Lemma \ref{lemmacs} that $u(\nabla, \nabla_0, h)$ are differential forms. From Stokes formula
(\ref{Stokes}) it immediately follows that they are closed. To prove that they are real
we use (iii) of the previous Lemma:
\begin{eqnarray}  
   & \overline{u(\nabla, \nabla_0, h)}= 
     (-i)^{p+1}(\overline{cs_{p}(\nabla, \nabla_0)}+ \overline{cs_{p}(\nabla_{0}, \nabla_{0}^{h})}+ 
     \overline{cs_{p}(\nabla_{0}^{h}, \nabla^{h})})  = & \nonumber \\
   & + (-i)^{p+1}(-1)^{p} cs_{p}(\nabla^{h}, \nabla_{0}^{h})+ cs_{p}(\nabla_{0}^{h}, \nabla_{0})+ 
     cs_{p}(\nabla_{0}, \nabla)) =  & \nonumber \\
   & =  (-i)^{p+1}(-1)^{p}(-1) u(\nabla, \nabla_0, h)= u(\nabla, \nabla_0, h) & \nonumber 
\end{eqnarray}
\hspace*{.3in}If $\nabla_1$ is another connection,
using  (\ref{Stokes}) again, it follows that $u(\nabla, \nabla_0, h)- u(\nabla, \nabla_1, h)$ $=$ $i^{p+1}dv$ where $v$ is the (linear!) differential form 
\[ v= cs_{p}(\nabla, \nabla_0, \nabla_1)- cs_{p}(\nabla^h, \nabla^{h}_{0}, \nabla^{h}_{1})+ cs_{p}(\nabla_0, \nabla_{0}^{h}, \nabla_1)-
cs_{p}(\nabla_{0}^{h}, \nabla_1, \nabla_{1}^{h}) \ .\]
(iii) clearly follows from (ii), which in turn follows from (ii) of Lemma \ref{lemmacs} and the fact that $\nabla\sim \nabla_{0}$ implies $\nabla^{h}\sim \nabla_{0}^{h}$. To see that our classes do not depend on $h$, it suffices to show that given a linear connection 
$\nabla$ on a vector bundle $F$, $cs_{p}(\nabla, \nabla^{h})$ is independent of $h$ up to the boundary of a differential form.
Let $h_0$ and $h_1$ be two metrics. Although the proof below works for general $\nabla$'s, simpler formulas are possible when $\nabla$ is flat. So, let us first assume that (actually we may assume that $\nabla$ is the canonical connection on a trivial vector bundle). From Stokes' formula (\ref{Stokes}) applied to $(\nabla, \nabla^{h_0}, \nabla^{h_1})$, it suffices to
show that $cs_{p}(\nabla^{h_0}, \nabla^{h_1})$ is a closed form. We choose a family $h_{t}$ of metrics joining $h_{0}$ and $h_{1}$. 
One only has to show
that $\frac{\partial}{\partial t} cs_{p}(\nabla^{h_0}, \nabla^{h_t})$ are closed forms. Writing $h_{t}(x, y)= h_{0}(u_{t}(x), y)$,
these Chern-Simons forms are, up to a constant, $Tr(\omega_{t}^{2p-1})$ 
where 
\[ \omega_{t}= \nabla^{h_t}- \nabla^{h_0} = u_{t}^{-1} d_{\nabla^{h_0}}(u_t) \]
(here is where we use the flatness of $\nabla$). A simple computation shows that 
\[ \frac{\partial \omega_t}{\partial t}= d_{\nabla^{h_0}}(v_t)+ [\omega_t, v_t]\ ,\]
where $v_t= u_{t}^{-1}\frac{\partial u_{t} }{\partial t}$. Since $d_{\nabla^{h_0}}(\omega_{t}^{2})= 0$, this implies
\[ \frac{\partial \omega_t}{\partial t}\omega_{t}^{2p-2}= d_{\nabla^{h_0}}(v_t\omega_{t}^{2p-2})+ [\omega_t, v_t\omega_{t}^{2p-2}] \ .\]
Now, by the properties of the trace it follows that
\[ \frac{\partial}{\partial t} Tr_{s}(\omega_{t}^{2p-1})= d Tr_{s}(v_{t}\omega_{t}^{2p-2})  \]
as desired. Assume now that $\nabla$ is not flat. We choose a vector bundle $F{\,'}$ together with a connection $\nabla^{\,'}$
compatible with a metric $h^{\,'}$, such that $\tilde{F}= F\oplus F^{\,'}$ admits a flat connection $\nabla_0$. 
We put $\tilde{\nabla}= \nabla\oplus \nabla^{\,'}$ and,
for any metric $h$ on $F$, we consider the metric $\tilde{h}= h\oplus h^{\,'}$ on $\tilde{F}$. Clearly
$cs_{p}(\tilde{\nabla}, \tilde{\nabla}^{\tilde{h}})= cs_{p}(\nabla, \nabla^{h})$. Using also (iii) of Lemma \ref{lemmacs} and
Stokes' formula, we have:
\begin{eqnarray}
 & cs_{p}(\nabla, \nabla^{h}) & = cs_{p}(\nabla_{0}, \nabla_{0}^{\tilde{h}})- cs_{p}(\nabla_{0}, \tilde{\nabla})+ (-1)^{p}\overline{cs_{p}(\nabla_{0}, \tilde{\nabla})}  \nonumber \\
 &  & + d(cs_{p}(\nabla_0, \tilde{\nabla}, \tilde{\nabla}^{\tilde{h}})- cs_{p}(\nabla_{0}, \tilde{\nabla}_{0}, \tilde{\nabla}^{\tilde{h}}))  .
\nonumber
\end{eqnarray} 
Hence, by the flat case, $cs_{p}(\nabla, \nabla^{h})$ modulo exact forms does not depend on $h$.\\
For (iv) one uses Stokes' formula (\ref{Stokes}) and (ii) of Lemma \ref{lemmacs} to conclude
that $cs_{p}(\nabla^{\,'}, \nabla_{0})- cs_{p}(\nabla, \nabla_{0})$ is the differential of the linear form $cs_{p}(\nabla, \nabla^{\,'}, \nabla_{0})$. To prove (v) we only have to show (see (i)) that there exists a linear connection $\nabla^{0}$
on $E$ which is compatible with both $\partial$ and $\partial^{\,'}$. For this, one defines $\nabla^{0}$ locally by $\nabla^{0}_{f \frac{\partial}{\partial x_{k}}}= f\nabla_{\frac{\partial}{\partial x_{k}}}$, and then use a partition of unity argument.\\
We now assume that $E$ is real. From Lemma \ref{lemmacs}, 
\[ cs_{p}(\nabla_{m}, \nabla_{0}^{h})= (-1)^{p} cs_{p}(\nabla_{m}^{h}, \nabla_{0})= (-1)^{p+1}cs_{p}(\nabla_{0}, \nabla_{m})\ .\]
Combined with Stokes' formula (\ref{Stokes}), this implies 
\[ d cs_{p}(\nabla_{0}, \nabla_{m}, \nabla_{0}^{h})= (1+ (-1)^{p+1}) cs_{p}(\nabla_{0}, \nabla_{m})- cs_{p}(\nabla_{0}, \nabla_{0}^{h})\ ,\]
which proves (vi).\ \ \ $\Boxe$\\

\hspace*{.1in}Note that the construction of the flat characteristic classes presented here actually works for $\nabla$'s which are 
``flat up to homotopy'', i.e. whose curvatures are of type 
$[-, \partial ]$. Moreover, this notion is stable under equivalence, and the flat characteristic classes only depend on the equivalence
class of $\nabla$ (cf. (iv) of the Theorem). Note also that, as in \cite{Crai} (and following \cite{BiLo}), there is a version of our discussion for 
super-connections \cite{Chern} up to homotopy. Some of our constructions can then be interpreted in terms of the super-connection $\partial+ \nabla$.\\
\hspace*{.3in}If $E$ is regular in the sense that $Ker(\partial)$ and $Im(\partial)$ are vector bundles,
then so is the cohomology $H(E, \partial)= Ker(\partial)/Im(\partial)$, and any connection up to homotopy
$\nabla$ on $(E, \partial)$ defines a linear connection $H(\nabla)$ on $H(E)$. Moreover,
$H(\nabla)$ is flat if $\nabla$ is, and the characteristic classes $u_{2p-1}(E, \partial, \nabla)$ probably coincide with
the classical \cite{BiLo, KaTo} characteristic classes of the flat vector bundle
$H(E, \partial)$. 
In general, the $u_{2p-1}(E, \partial, \nabla)$'s should be viewed as invariants of $H(E, \partial)$
constructed in such a way that no regularity assumption is required. Let us discuss here an instance 
of this. We say that $E$ is $\mathbb{Z}$-graded if it comes from a cochain complex
\begin{equation}\label{lasteq} 
0\rmap E(0)\stackrel{\partial}{\rmap} E(1) \stackrel{\partial}{\rmap} \ldots \stackrel{\partial}{\rmap} E(n) \rmap 0\ ,
\end{equation}
In other words, it must be of type  
$E= \oplus_{k=0}^{n}E(k)$ with the even/odd $\mathbb{Z}_{2}$-grading, and with $\partial(E(k))\subset E(k+1)$. As usual,
we say that $E$ is acyclic if $Ker(\nabla)= Im(\nabla)$ (i.e. if (\ref{lasteq}) is exact).

\begin{prop}\label{propr2} $\text{}$ 
\begin{enumerate}[(i)]
\item If $(E, \partial)$ is acyclic, then any two connections up to homotopy on $(E, \partial)$ are equivalent.
Moreover, if $E$ is $\mathbb{Z}$-graded, then $u_{2p-1}(E, \partial, \nabla)= 0$.
\item If $(E^{k}, \partial^{k}, \nabla^{k})$ are $\mathbb{Z}$-graded complexes of vector bundles endowed with flat connections up to homotopy
which fit into an exact sequence
\begin{equation}\label{theseq}
0\rmap E^0\stackrel{\delta}{\rmap} E^{1}\stackrel{\delta}{\rmap}\ldots \stackrel{\delta}{\rmap} E^{n}\rmap 0 
\end{equation}
compatible with the structures (i.e. $[\delta, \partial ]= [\delta, \nabla ]= [\delta, H_{\nabla}]= 0$), then
\[ \sum_{k=0}^{n} (-1)^{k} u_{2p-1}(E^{k}, \partial^{k}, \nabla^{k})= 0\ .\]
\end{enumerate}
\end{prop}

{\it Proof:} The second part follows from (i) above and (v) of Theorem \ref{theorem2}. To see this, we form the 
super-vector bundle $E= \oplus_{k} E^{k}$ (which is $\mathbb{Z}$-graded by the total degree) and the direct sum (non-linear) connection $\nabla$ acting on $E$. Then
$\nabla$ is a connection up to homotopy in both $(E, \partial)$ and $(E, \partial+ \delta)$. Clearly
$u_{2p-1}(E, \partial, \nabla)=\sum_{k=0}^{n} (-1)^{k} u_{2p-1}(E^{k}, \partial^{k}, \nabla^{k})$,
while the exactness of (\ref{theseq}) implies that $\partial+ \delta$ is acyclic.
Hence we are left with (i). For the first part we remark that the acyclicity assumption implies that $\partial^*\partial+ \partial\partial^*$ 
is an isomorphism (``Hodge''). Then any operator $u$ which commutes with $\partial$ can be written as a commutator
$[-v, \partial]$ where 
\begin{equation}\label{Hodge}
v= ua, \ \ a= -(\partial^*\partial+ \partial\partial^*)^{-1}\partial^*\ .
\end{equation}
This applies in particular to $u= \nabla^{\,'}-\nabla$
for any two connections up to homotopy on $(E, \partial)$. We now have to prove that $cs_{p}(\nabla, \nabla^{h})$ is
zero in cohomology, where $\nabla$ is a linear connection on $(E, \partial)$, and $h$ is a metric. 
For this we use a result of \cite{BiLo} (Theorem 2.16) which says that $cs_{p}(A, A^{h})$ are closed forms 
provided $A= A_{0}+ A_{1}+ A_{2} + \ldots$ is a flat super-connection \cite{Chern} on $E$ with the properties:
\begin{enumerate} [(i)]
\item $A_1$ is a connection on $E$ preserving the $\mathbb{Z}$-grading,
\item $A_k\in \A^{k}(M; \text{Hom}(E^{*}, E^{*+1-k}))$ for $k\neq 1$.
\end{enumerate}

\begin{lem}\label{Fedosov} Given a (linear) connection $\nabla$ on the acyclic cochain complex (\ref{lasteq}),
there exists a  super-connection on $E$ of type
\[ A= \partial + \nabla + A_{2}+ A_{3}+ \ldots \ \, : \A(M; E)\rmap \A(M; E) \ ,\]
which is flat and satisfies (i) and (ii) above.
\end{lem}

Let us show that this lemma, combined with the result of \cite{BiLo} mentioned above, prove the desired result.
Using Stokes' formula it follows that 
\begin{eqnarray}
 & cs_{p}(\nabla, \nabla^h)& = cs(A, A^{h})+ d( cs_{p}(\nabla, \nabla^h, A^h)- cs_{p}(\nabla, A, A^{h}))+\nonumber \\
 &    & + cs_{p}(\nabla, A)- cs_{p}(\nabla^h, A^h) \ ,\nonumber
\end{eqnarray}
and we show that $cs_{p}(\nabla, A)= 0$ (and similarly that $cs_{p}(\nabla^h, A^h)= 0$). 
Writing $\theta= A- \nabla$ and using the definition of the Chern-Simons forms, it suffices to prove that
\[ Tr_{s}(((1-t^2)\nabla^2+ (t- t^2)[\nabla, \theta])^{p-1}\theta)= 0\]
for any $t$. Since the only endomorphisms of $E$ which count are those preserving the degree, we see that the only 
term which can contribute is $Tr_{s}(\nabla^{2(p-2)}[\nabla, \theta]\theta)$ $=$ $Tr_{s}(\nabla^{2(p-2)}[\nabla, A_{2}]\partial)$. 
But $\nabla^{2(p-2)}[\nabla, A_{2}]\partial$ commutes with $\partial$ hence its super-trace must vanish (since $Tr_s$ commutes with taking cohomology). \ \ $\Boxe$\\

{\it Proof of Lemma \ref{Fedosov}:} (Compare with \cite{Fed}). The flatness of $A$ gives us certain equations on the $A_k$'s that we can solve inductively,
using the same trick as in (\ref{Hodge}) above. For instance, the first equation
is $[\partial, A_{2}]+ \nabla^{2}= 0$. Since $u_{1}= \nabla^{2}$ commutes with $\partial$, this equation
will have the solution $A_{2}= u_{1}a$ (with $a$ as in (\ref{Hodge})). The next equation is $[\partial, A_{3}]+
[A_{1}, A_{2}]= 0$. It is not difficult to see that $u_{2}= [A_{1}, A_{2}]$ commutes with $\partial$, and we
put $A_{3}= u_{2}a$. Continuing this process, at the $n$-th level we put $A_{n+1}= u_{n}a$
where $u_{n}= [\nabla, A_{n}]+ [A_{1}, A_{n-2}]+ \ldots $ as they arise from the coresponding equation. We
leave to the reader to show that the $u_{n}$'s also satisfy the equations
\[ u_{n}= u_{n-1}[\nabla, a]+ (\sum_{i+j= n-1} u_{i}u_{j})a^2  .\]
Since $[\partial, a]= -1$, $\partial$ will commute with both $[\nabla, a]$ and $a^2$, hence also with the $u_{n}$'s (induction on $n$).  
It then follows that $A_{n+1}$ satisfies the desired equation $[\partial, A_{n+1}]= -u_n$. \ \ $\Boxe$\\


\section*{Application to Lie algebroids} 

Recall \cite{McK} that a {\it Lie algebroid} over $M$ consists of a Lie bracket $[\cdot \, , \cdot ]$
defined on the space $\Gamma \mathfrak{g}$ of sections of a vector bundle $\mathfrak{g}$ over $M$,
together with a morphism of vector bundles $\rho: \mathfrak{g}\rmap TM$ ({\it the anchor}
of $\mathfrak{g}$) satisfying $[X, fY] = f[X,Y]+ \rho(X)(f) \cdot Y$ for all
$X, Y \in \Gamma(\mathfrak{g})$ and $f \in C^{\infty}(M)$. Important examples are 
tangent bundles, Lie algebras, foliations, and algebroids associated to Poisson manifolds.
It is easy to see (and has already been remarked in many other places \cite{McK}, \cite{Cra}, \cite{Fer1}, etc. etc.)
that many of the basic constructions involving vector fields have a straightforward $\mathfrak{g}$-version
(just replace $\xx(M)$ by $\Gamma(\mathfrak{g})$).
Let us briefly point out some of them.
\begin{enumerate} [(a)]
\item {\it Cohomology:} the Lie-type formula (\ref{differential}) for the classical De Rham differential makes sense
for $X\in \Gamma\mathfrak{g}$ and defines a differential $d$ on the space $C^{*}(\mathfrak{g})= \Gamma \Lambda^*\mathfrak{g}^*$,
hence a cohomology theory $H^{*}(\mathfrak{g})$. Particular cases are De Rham cohomology, Lie algebra cohomology, foliated
cohomology, and Poisson cohomology.
\item {\it Connections and Chern characters:} According to the general philosophy, $\mathfrak{g}$-connections
on a vector bundle $E$ over $M$ are linear maps $\Gamma(\mathfrak{g})\times \Gamma E\rmap \Gamma E$ satisfying the usual 
identities. Using their curvatures, one obtains $\mathfrak{g}$- Chern classes $Ch^{\mathfrak{g}}(E)\in H^*(\mathfrak{g})$
independent of the connection.
\item {\it Representations:} Motivated by the case of Lie algebras, and also by the relation to groupoids (see e.g. \cite{Cra}),
vector bundles $E$ over $M$ together with a flat $\mathfrak{g}$-connection are called representations of $\mathfrak{g}$.
This time $\nabla$ should be viewed as an (infinitesimal) action of $\mathfrak{g}$ on $E$. 
\item {\it Flat characteristic classes:} The explicit approach to flat characteristic classes 
(as e.g. in \cite{BiLo}, or as in the previous section) has an obvious  $\mathfrak{g}$-version.
Hence, if $E$ is a representation of $\mathfrak{g}$, then $Ch^{\mathfrak{g}}(E)= 0$,
and one obtains the secondary characteristic classes $u_{2p-1}(E)\in H^{2p-1}(\mathfrak{g})$.
Maybe less obvious is the fact that one can also extend the Chern-Weil type approach, at the level of frame bundles
(as e.g. in \cite{KaTo}). This has been explained in \cite{Cra}, and has certain advantages 
(e.g. for proving ``Morita invariance'' of the $u_{2p-1}(E)$'s and for relating them to differentiable cohomology). 
\item {\it Up to homotopy:} All the constructions and results of the previous sections carry over to Lie algebroids without any problem.
As above, a {\it representation up to homotopy} of $\mathfrak{g}$ is a supercomplex (\ref{complex}) of vector bundles
over $M$, together with a flat $\mathfrak{g}$-connection up to homotopy. 
\item {\it The adjoint representation:} The main reason for working ``up to homotopy'' is that
the adjoint representation of $\mathfrak{g}$ only makes sense as a representation up to homotopy \cite{ELW}. Roughly
speaking, it is the formal difference $\mathfrak{g}- TM$. The precise definition is:
\begin{equation}\label{adjrep}
\text{Ad}(\mathfrak{g}): \ \ \ \xymatrix{  \mathfrak{g}\ \ar@<-1ex>[r]_-{\rho} & \ \ TM \ar@<-1ex>[l]_-{0}\ \ \ \ , }
\end{equation}
with the flat $\mathfrak{g}$-connection up to homotopy $\nabla^{ad}$ given by 
$\nabla_{X}^{ad}(Y)= [X, Y]$, $\nabla_{X}^{ad}(V)= [\rho(X), Y]$
(and the homotopies $H(f, X)(Y)= 0$, $H(f, X)(V)= V(f) X$), for all $X, Y\in \Gamma\mathfrak{g}$, $V\in \xx(M)$.
\end{enumerate}

\hspace*{.1in}Let us denote by $u_{2p-1}^{\mathfrak{g}}$ the characteristic classes $u_{2p-1}(\text{Ad}(\mathfrak{g}))$ of the adjoint representation. 
The most useful description from a computational (but not conceptual) point of view is given by (vi) of Theorem \ref{theorem2}
(more precisely, its $\mathfrak{g}$-version). 

\begin{num}{\bf Definition}\rm\ We call {\it basic $\mathfrak{g}$-connection} 
any $\mathfrak{g}$-connection on $\text{Ad}(\mathfrak{g})$ which is equivalent to the adjoint connection $\nabla^{\text{ad}}$. 
\end{num}

It is not difficult to see that any such connection is also basic in sense of \cite{Fer1}
(and the two notions are equivalent at least in the regular case). Hence we have the following possible description of the $u_{2p-1}^{\mathfrak{g}}$'s,
which shows the compatibility with Fernandes' intrinsic characteristic classes \cite{Fer1, Fer2}:
\[ u_{2p-1}^{\mathfrak{g}}= \left \{ \begin{array}{ll}
                                               0 & \mbox{if $p=$ even} \\
                                               \frac{1}{2}(-1)^{\frac{p+1}{2}} cs_{p}(\nabla_{\text{bas}}, \nabla_{\text{m}}) & \mbox{if $p=$ odd}
                                                        \end{array}
                                            \right. ,\]
where $\nabla_{\text{bas}}$ is any basic $\mathfrak{g}$-connection, and $\nabla_{\text{m}}$ is any metric connection on $\mathfrak{g}\oplus TM$.
Hence the conclusion of our discussion is the following (which can also be taken as definition of the characteristic classes of \cite{Fer1, Fer2}).

\begin{st} If $E$ is a representation up to homotopy then 
$Ch^{\mathfrak{g}}(E)= 0$, and the secondary characteristic classes $u_{2p-1}(E)\in H^{2p-1}(\mathfrak{g})$ of representations \cite{Crai}
can be extended to such representations up to homotopy. When applied to the adjoint representation ${\rm Ad}(\mathfrak{g})$, the resulting
classes $u_{2p-1}^{\mathfrak{g}}$ are (up to a constant) the intrinsic characteristic classes of $\mathfrak{g}$ \cite{Fer1}.
\end{st}

\hspace*{-.2in}{\bf More on basic connections:} Let us try to shed some light on the notion of basic $\mathfrak{g}$-connection. In our context these are the linear connections which are equivalent to the adjoint connection, while in \cite{Fer1} they appear as a natural extension of Bott's basic connections for foliations. Although not flat in general, they are always flat up to homotopy. The existence
of such connections is ensured by Lemma \ref{lemmaexist} and it was also proven in \cite{Fer1}. There is however
a very simple and explicit way to produce them out of ordinary connections on the vector bundle $\mathfrak{g}$.

\begin{prop} Let $\nabla$ be a connection on the vector bundle $\mathfrak{g}$. Then the formulas
\begin{eqnarray}
&  \check{\nabla}^{0}_{X}(Y)= [X, Y]+ \nabla_{\rho(Y)}(X) \nonumber \\
&  \check{\nabla}^{1}_{X}(V)= [\rho(X), V]+ \rho(\nabla_{V}(X)) \nonumber
\end{eqnarray}
($X, Y\in \Gamma \mathfrak{g}$, $V\in \Gamma TM$) define a basic $\mathfrak{g}$-connection
$\check{\nabla}= (\check{\nabla}^0, \check{\nabla}^1)$.
\end{prop}

{\it Proof:} We have $\check{\nabla}= \nabla^{ad}+ [\theta, \partial ]$, where $\theta$ is the (non-linear)
$\text{End}(\text{Ad}(\mathfrak{g}))$-valued form on $\mathfrak{g}$ given by $\theta(X)(V)= \nabla_{V}(X)$, $\theta(X)(Y)= 0$. \ \ $\Boxe$\\

\hspace*{.1in}Depending on the special properties of $\mathfrak{g}$, there are various other useful basic connections.
This happens for instance when $\mathfrak{g}$ is regular, i.e. when the rank of the anchor $\rho$ is constant. 
Let us argue that, in this case, the adjoint representation is (up to homotopy) the formal difference $K- \nu$,
where $K$ is the kernel of $\rho$, and $\nu$ is the normal bundle $TM/\F$ of the foliation $\F= \rho(\mathfrak{g})$.
This time, Bott's formulas \cite{Bott} trully make sense on $\nu$ and $K$, making them into honest representations of $\mathfrak{g}$:
\begin{eqnarray}
& \nabla_{X}(\bar{Y})= \overline{[X, Y]}, \ \ \ \forall\ X\in \Gamma\mathfrak{g},\ \bar{Y}\in\Gamma\nu \label{Bottrep1} \\
& \nabla_{X}(Y)= [X, Y], \ \ \ \forall\ X\in \Gamma\mathfrak{g},\  Y\in \Gamma K\ \ . \label{Bottrep2}
\end{eqnarray}
Now, choosing splittings $\alpha: \F\rmap \mathfrak{g}$ for $\rho$, and $\beta: TM\rmap \F$ for the inclusion,
we have induced decompositions
\[ \mathfrak{g}\cong  K\oplus \F, \ \ \ TM\cong \nu\oplus \F\ .\]
As mentioned above, the formal difference $K- \nu$ (view it as a graded complex with $K$ in even degree, $\nu$ in odd degree,
and zero differential) is a representation of $\mathfrak{g}$. On the other hand, any $\F$-connection $\nabla$ on $\F$
defines a $\mathfrak{g}$-connection on the super-complex
\[ D(\F): \xymatrix{  \F\ \ar@<-1ex>[r]_-{id} & \ \ \F \ar@<-1ex>[l]_-{0} }\]
(and its homotopy class does not depend on $\nabla$). Hence one has an induced $\mathfrak{g}$-connection $\nabla^{\alpha, \beta}$ on $\text{Ad}(\mathfrak{g})$, so that $(\text{Ad}(\mathfrak{g}), \nabla^{\alpha, \beta})$ is isomorphic to $(K- \nu) \oplus D(\F)$.
Explicitly,
\begin{eqnarray}
&  \nabla^{\alpha, \beta}_{X}(Y)= [X, Y- \alpha\rho(Y)]+ \alpha\nabla_{\rho(Y)}(\rho{X}) \nonumber \\
&  \nabla^{\alpha, \beta}_{X}(V)= [\rho(X), V]- \beta [\rho(X), V]+ \nabla_{\rho(X)}(\beta(V)) \nonumber
\end{eqnarray}
for all $X, Y\in \Gamma\mathfrak{g}$, $V\in \xx(M)$. Note that the second part of the following proposition can also be derived
from (iv) of Proposition \ref{propr2}.

\begin{prop} Assume that $\mathfrak{g}$ is regular. For any $\F$-connection $\nabla$ on $\F$, and any splittings $\alpha$, $\beta$ as above,
$\nabla^{\alpha, \beta}$ is a basic $\mathfrak{g}$-connection. In particular 
\[ u_{2p-1}^{\mathfrak{g}}= u_{2p-1}(K)- u_{2p-1}(\nu) \ ,\]
where $K$ and $\nu$ are the representations of $\mathfrak{g}$ defined by Bott's formulas (\ref{Bottrep1}), (\ref{Bottrep2}).
\end{prop}

{\it Proof:} We have $\nabla^{\alpha, \beta}= \nabla^{\text{ad}}+ [\theta, \partial ]$, where $\theta$ is the $\text{End}(\text{Ad}(\mathfrak{g}))$-valued
non-linear form which is given by
\[ \theta(X)(V)= \alpha[\rho(X), \beta(V)]- \alpha\beta [\rho(X), V]- [X, \alpha\beta(V)]+ \alpha\nabla_{\rho(X)}\beta(V) \]
for $V\in \Gamma(TM)$, while $\theta(X)= 0$ on $\mathfrak{g}$ (we leave to the reader to check that the previous formula
is $C^{\infty}(M)$-linear on $V$). \ \ $\Boxe$.\\

Marius Crainic,\\
\hspace*{.2in}Utrecht University, Department of Mathematics,\\
\hspace*{.2in}P.O.Box:80.010,3508 TA Utrecht, The Netherlands,\\ 
\hspace*{.2in}e-mail: crainic@math.ruu.nl\\
\hspace*{.2in}home-page: http://www.math.uu.nl/people/crainic/

\end{document}